\documentclass[leqno,12pt]{article}
\usepackage{amssymb,amsfonts}
\usepackage{amsmath,latexsym}

\newcommand{\bepr}{{\em Proof} } 
\newcommand{\enpr}{\hfill \rule{.5em}{.5em}}

\newcommand{\R}{{\mathbb R}}

\newcommand{\C}{{\mathbb C}}

\newtheorem{defin}{Definition}[section] 
\newtheorem{prop}{Proposition}[section] 
\newtheorem{thm}{Theorem}[section] 
\newtheorem{lemma}{Lemma}[section] 

\newtheorem{cor}{Corollary}[section]

\begin{document}

\title{Expansion of a compressible gas in vacuum}

\author{Denis Serre\thanks{UMPA, UMR 5669 CNRS-ENS Lyon~;
46, all\'ee d'Italie 
69364 LYON Cedex 07, FRANCE}}

\date{\today}

\maketitle

\centerline{\em Dedicated to Tai-Ping Liu, on the occasion of his 70th birthday}

\centerline{\em  En m\'emoire de G\'erard Lasseur}

\begin{abstract}
Tai-Ping Liu \cite{Liu_JJ} introduced the notion of ``physical solution'' of the isentropic Euler system when the gas is surrounded by vacuum. This notion can be interpreted by saying that the front is driven by a force resulting from a H\"older singularity of the sound speed. We address the question of when this acceleration appears or when the front just move at constant velocity.

We know from \cite{Gra,SerAIF} that smooth isentropic flows with a non-accelerated front exist globally in time, for suitable initial data. In even space dimension, these solutions may persist for all $t\in\R$~; we say that they are {\em eternal}. We derive a sufficient condition in terms of the initial data, under which the boundary singularity must appear. As a consequence, we show  that, in contrast to the even-dimensional case, eternal flows with a non-accelerated front don't exist in odd space dimension.

In one space dimension, we give a refined definition of physical solutions. We show that for a shock-free flow, their asymptotics as both ends $t\rightarrow\pm\infty$ are intimately related to each other.
\end{abstract}

\section{The isentropic Euler equations}

We consider a compressible isentropic gas of density $\rho$, velocity $u$, pressure $p$ and specific internal energy $e$. The flow is governed by the Euler system
\begin{eqnarray}
\label{eq:cmass} \partial_t\rho+{\rm div}(\rho u) & = & 0, \\
\label{eq:cmom} \partial_t(\rho u)+{\rm Div}(\rho u\otimes u)+\nabla p(\rho) & = & 0.
\end{eqnarray}
Classical solutions of (\ref{eq:cmass},\ref{eq:cmom}) satisfy in addition the conservation of energy
$$\partial_t(\frac12\rho |u|^2+\rho e)+{\rm div}((\frac12\rho |u|^2+\rho e+p) u) =0.$$
When the solution develops singularities, say shock waves, the latter equality is replaced by the inequality
\begin{equation}\label{eq:ener}
\partial_t(\frac12\rho |u|^2+\rho e)+{\rm div}((\frac12\rho |u|^2+\rho e+p) u) \le0
\end{equation}
which holds true in the distributional sense~; (\ref{eq:ener}) plays the role of an admissibility condition.

For the sake of clarity, we assume a polytropic equation of state
$$p(\rho)=A\rho^\gamma,\qquad A>0,$$
where $\gamma>1$ is the {\em adiabatic constant}. Then the internal energy is given by
\begin{equation}\label{eq:rhoe}
\rho e=\frac p{\gamma-1}\,.
\end{equation}
The sound speed is  $c=\sqrt{p'}=\sqrt{\gamma A}\,\rho^\kappa$ where $\kappa:=\frac{\gamma-1}2\,$. We recall that a {\em mono-atomic} gas corresponds to the choice
$$\gamma_d:=1+\frac2d\,,$$
which means that the molecules have only $d$ degrees of freedom, all of them associated with translations in space. In the physically relevant case $d=3$, $\gamma_3=\frac53$ is the parameter associated with mono-atomic gases like Argon. Di-atomic gases, like $H_2,O_2$ or the air, obey to a pressure law with $\gamma=\frac75<\gamma_3$. If $d=2$, then $\gamma_2=2$ corresponds to the shallow water equations.

\bigskip

We are interested in this paper in flows for which the total mass and energy are finite, and the gas occupies a bounded region $\Omega(t)$ of the ambient space $\R^d$ at time $t$. The gas is surrounded by vacuum. We assume {\em a priori} that  the front $\Gamma(t)=\partial\Omega(t)$ is a smooth or piecewise smooth hypersurface in $\R^d$. We  denote
$$\Omega=\{(x,t)\,|\,x\in\Omega(t)\hbox{ and }t\in\R\},\qquad\Gamma=\{(x,t)\,|\,x\in\Gamma(t)\hbox{ and }t\in\R\}.$$
For weak admissible solutions of the Euler system, the total mass is conserved and the total energy is a non-increasing function of time~:
$$\int_{\R^d}\rho(x,t)\,dx\equiv\int_{\R^d}\rho_0(x)\,dx=: M<\infty,\qquad\frac{d}{dt}\int_{\R^d}(\frac12\rho |u|^2+\rho e)(x,t)\,dx\le0.$$
The latter inequality is an equality if the solution is classical. 

Of particular importance is the nature of the boundary condition along $\Gamma(t)$.
On the one hand, because there should not be any transfer of mass accross the boundary, $\Gamma$  moves at the velocity $u$. More precisely, its normal velocity equals $u\cdot \nu$, where $\nu$ is the unit outer normal. The boundary can be viewed as a collection of particles moving at the fluid velocity $u$~; their paths are integrals curves of the ODE
\begin{equation}\label{eq:ode}
\frac{dX}{dt}=u(X(t),t).
\end{equation}
On the other hand, the Rankine--Hugoniot relations reduce on $\Gamma$ to $[p(\rho)]=0$, where
$\rho$ vanishes on the vacuum side~; therefore $\rho$ must vanish on the interior side, meaning that $\rho$ is continuous across $\Gamma(t)$ (see also \cite{LS}). We mention in passing that this property might be violated at isolated times because an isolated (in space {\em and} time) discontinuity is not a shock front. An explicit example of that possibility was given by Greenspan \& Butler \cite{GB}, see Section \ref{s:phy} (``sub/supersonic'' example). More generally, we allow discontinuities that occur on a codimension-$1$ subset of $\Gamma$, because they are not shock waves.

Once we know that $\rho\equiv0$ along the boundary, there remains to understand its regularity, or singularity, in terms of the distance to the $\Gamma(t)$. Here one can think of two approaches, complementary to each other. Both are based on the quasi-linear symmetric form of the system,
\begin{eqnarray}
\label{eq:qlmass} \partial_t\bar c+(u\cdot\nabla)\bar c+\kappa\bar c\,{\rm div}u & = & 0, \\
\label{eq:qlmom} \partial_t u+(u\cdot\nabla)u+\kappa\bar c\nabla \bar c & = & 0,
\end{eqnarray}
where $\bar c:=\frac c\kappa\,$ is a renormalized version of the sound speed.  We point out that system (\ref{eq:qlmass},\ref{eq:qlmom}) removes the singularity at $\rho=0$, inherent to the Euler system, whose hyperbolicity degenerates at vacuum. These approaches distinguish two regimes, depending on whether  the quantity $\bar c\nabla\bar c$ vanishes or not. If it does, then the particles at the front move freely~; we say that the front is not accelerated, or that the flow is {\em smooth up to vacuum}. If instead $c\nabla c\ne0$, then the front experiences a normal acceleration. This is what T.-P. Liu called a {\em physical singularity} at the boundary. Of course it requires that $c$ be $\frac12$-H\"older at the boundary. Beware that a given flow can be smooth up to vacuum in some region, for instance on some time interval, and display a physical singularity elsewhere. The transition between both regimes is still not completely understood.


\paragraph{Outline of the paper.} Section \ref{s:smooth} deals with flows for which the sound speed is smooth enough at the vacuum, so that the front is not accelerated. Such flows arise for convenient initial data, but might exist only for a finite time. When the boundary regularity is lost, it is expected that a H\"older-type singularity of $c$ develops at the boundary. These singularities, called {\em physical} by T.-P. Liu, are investigated in Section \ref{s:phy}.

Our main results are Theorem \ref{th:nonex} (non-existence in odd space dimension of eternal flows that are smooth up to vacuum) and Theorem \ref{th:size} on the relation between the asymptotics as $t\rightarrow\pm\infty$ in one space-dimension.

\section{Flows that are smooth up to vacuum}\label{s:vacuum}

The first approach considers  a solution $(c,u)$ of (\ref{eq:qlmass},\ref{eq:qlmom}), smooth   in the entire space $\R^d$, with a compactly supported component $c$. The following result is a  classical application of the theory of symmetric hyperbolic systems, which can be found in \cite{Daf}.
\begin{thm}\label{th:smooth}
Let the initial data $(c_0,u_0)$ belong to the Sobolev space $H^s(\R^d)$ for some parameter $s>1+\frac d2\,$. Then there exists an open time interval $I$ containing $t=0$, and a unique local-in-time solution in the class ${\cal C}(I;H^s)\cap{\cal C}^1(I;H^{s-1})$. 
\end{thm}
Thanks to Sobolev embedding, this solution is classical~: $(c,u)\in{\cal C}^1(I\times\R^d)$. If $c_0$ is compactly supported, then $c(\cdot,t)$ is too, and $c$ vanishes at the boundary of its support.

\bigskip

The following argument is adapted from Liu \& Yang \cite{LY1}, where it was designed in the context of damped flows. For a solution given by Theorem \ref{th:smooth}, the equation (\ref{eq:qlmom}) reduces to
\begin{equation}\label{eq:trans}
\partial_t u+(u\cdot\nabla)u=0\qquad\hbox{along}\quad\Gamma. 
\end{equation}
Since $\Gamma$ is transported by the flow, this exactly means that the trajectories defined by (\ref{eq:ode})
and originating on the boundary (that is, $X(0)\in\Gamma(0)$) remain on $\Gamma$ and have a constant velocity
$$u(X(t),t)\equiv{\rm cst}.$$
We summarize this argument into
\begin{prop}\label{p:balist}
Let $(c,u)$ be a flow given by Theorem \ref{th:smooth} with $\Omega(t)$ bounded. Then the front $\Gamma(t)$ is obtained from $\Gamma (0)$ by a (uniform in time) transport:
$$\Gamma(t)=\psi_t(\Gamma(0))$$
where
$$\psi_t(x)=x+tu_0(x).$$
\end{prop}
We point out that $\psi$ is {\em not} the flow map of the gas, although it coincides with it on the vacuum boundary.

Because they are classical ones, the solutions provided by Theorem \ref{th:smooth} yield perfectly admissible solutions of the Euler equations (\ref{eq:cmass},\ref{eq:cmom}). At this stage, it seems that we don't need an extra condition at the boundary. Their only flaw is that they are defined only on a finite time interval, but this is something we are accustomate with in the theory of hyperbolic conservation laws. We anticipate that if $I=(-T_*,T^*)$ is bounded for a maximal solution, then some kind of singularity will appear as $t\rightarrow T^*$. We have in mind the formation of shock waves as usual, but because the front $\Gamma=\partial\Omega$ is a place where the hyperbolicity of the Euler system degenerates, the singularity might appear at $\Gamma(T^*)$ and then propagate along $\Gamma$. We shall see later on sufficient conditions that lead to such boundary singularities in finite time.

\subsection{Global smooth solutions}

It may happen that $T^*$ be infinite, in which case the smooth solution of the forward Cauchy problem is global-in-time.  Such solutions where first constructed in \cite{SerAIF} when $u_0$ is close to a linear field $x\mapsto Ax$, and the spectrum (the set of eigenvalues) of the matrix $A$ is contained in $\C\setminus(-\infty,0]$. This result was soon generalized to more realistic initial data~: 
\begin{thm}[M. Grassin \cite{Gra}.]\label{th:Mag}
In Theorem \ref{th:smooth}, assume in addition that the spectrum of $\nabla u_0(x)$ is contained in some fixed compact subset of $\C\setminus(-\infty,0]$. There exists an $\epsilon>0$, depending upon $u_0$, such that if 
$$\|c_0\|_{H^s}<\epsilon,$$
then the solution exists for all positive time~: $T^*=+\infty$.
\end{thm}

The idea behind the proof of Theorem \ref{th:Mag} is that the data is close to $(0,u_0)$. For the latter data, the solution is  $(c\equiv0,u)$, with $u$  governed by the vectorial Burgers equation
\begin{equation}\label{eq:Bur}
\partial_tu+(u\cdot\nabla)u=0.
\end{equation}
By assumption, the flow $\psi_t(x)=x+tu_0(x)$ is one-to-one at every positive time and therefore the solution of (\ref{eq:Bur}) is smooth for $t>0$. For a non-zero $c_0$, there is a competition between the dispersion induced by $\psi_t$ (the distance between two particles tends to increase linearly in time), and the nonlinear effect due to the pressure. Because of the dispersion, the density decays at an algebraic rate, and therefore the role of the pressure force gets smaller and smaller. If it was weak enough at initial time, we may expect that it will never be strong enough to lead to shock formation. An other important point is that the characteristic cones, which move at velocity $u+\xi$ with $|\xi|=c(\rho)$, will never overtake the boundary $\Gamma(t)$, because the slope of $c$ at the boundary is dominated by the diverging gradient of $u$.

\subsection*{Eternal smooth flows}

Theorem \ref{th:Mag}  has the beautiful consequence that smooth eternal compactly supported flows exist in even space dimension~:
\begin{cor}[$d$ even.]\label{c:even}
Suppose the space dimension $d$ is even.
Then in Theorem \ref{th:smooth}, one can choose the initial data so that the solution is eternal~: $I=\R$.
\end{cor}

\bepr

Just choose $u_0$ so that the spectrum of $\nabla u_0(x)$ is contained in a fixed compact subset $K$ of $\C\setminus\R$. Then choose $c_0$ small enough in $H^s(\R^d)$. Because $K$ does not meet $(-\infty,0]$, we have $T^*=+\infty$. Because $K$ does not meet $[0,+\infty)$, we have $T_*=-\infty$.

\enpr

\bigskip

Corollary \ref{c:even} raises the question of whether there exist eternal solutions with compact support, smooth up to vacuum, in odd space dimension. The construction in the proof above is not possible, because odd-size matrices do have at least one real eigenvalue. It is remarkable that one meets an obstruction, see Theorem \ref{th:nonex} below. The non-existence in one space dimension is a consequence of a calculus due to P. Lax \cite{Lax}~: the Euler system, whose wave velocities are $\lambda_\pm:=u\pm c$, can be diagonalized in Riemann coordinates
$$r_\pm=u\pm \bar c,$$
namely (\ref{eq:qlmass},\ref{eq:qlmom}) is equivalent to
$$(\partial_t+\lambda_\pm\partial_x)r_\pm=0.$$
Each of the transport equations above tells us that $r_\pm$ keeps a constant value along every characteristic curve $t\mapsto X(t)$ defined by (\ref{eq:ode}).
Lax showed that there exists a positive function $N(\rho)$ such that the expressions $y_\pm:=N(\rho)\partial_xr_\pm$ satisfy  Ricatti equations along the characteristics
\begin{equation}\label{eq:Ric}
(\partial_t+\lambda_\pm\partial_x)y_\pm+N^{-1}\frac{\partial\lambda_\pm}{\partial r_\pm}\,y_\pm^2=0.
\end{equation}
Because of 
genuine nonlinearity
$$\frac{\partial\lambda_\pm}{\partial r_\pm}>0,$$
the only eternal solution of (\ref{eq:Ric}) is $y_\pm\equiv0$. Therefore the global existence implies that $r_\pm$, or equivalently $\rho$ and $u$, are constant. With the assumption of finite mass, we conclude that $\rho\equiv0$~: there is no fluid at all. In the next section, we shall use the fact that if a solution $y$ of (\ref{eq:Ric}) exists on a time interval $(T,+\infty)$ (respectively on $(-\infty,T)$), then $y\ge0$, that is $\partial_xr_\pm\ge0$ (resp. $\partial_xr_\pm\le0$).

Notice that if $\gamma=\gamma_1=3$, the system above decouples as two independent Burgers equations, since then $r_\pm=\lambda_\pm$~:
\begin{equation}\label{eq:Burg}
(\partial_t+\lambda_\pm\partial_x)\lambda_\pm=0.
\end{equation}
In this case, one can take $N\equiv1$.

We summarize this analysis into 
\begin{prop}[$d$=1.]\label{p:und}
In one space dimension, there does not exist a non-trivial eternal solution $(c,u)$ of  (\ref{eq:qlmass},\ref{eq:qlmom}), of class ${\cal C}^2$, with $c$ of compact support.
\end{prop}

\subsection{Dispersion {\em vs} smooth boundary}\label{s:smooth}

The following statement tells us that if the initial velocity field does not drive fast enough the boundary, then a singularity must happen in finite time.
\begin{thm}\label{th:unilat}
We assume $\gamma\le1+\frac1{d-1}\,$.

Let $(c_0,u_0)\in H^s(\R^d)$ be an initial data, with  $s>1+\frac d2\,$. Assume that $c_0$ is compactly supported in $\Omega(0)$, a domain with smooth boundary. Assume that 
\begin{equation}\label{eq:jac}
J:=\int_{\Omega(0)}\det\nabla u_0(x)\,dx\le0.
\end{equation}
Then the maximal time $T^*>0$ of existence provided by Theorem \ref{th:smooth} is finite:
$$T^*<+\infty.$$
\end{thm}

Remark that the integral $J$ in (\ref{eq:jac}) depends only upon the restriction of $u_0$ to the boundary $\partial\Omega(0)$. The identity
$$\det\nabla u_0(x)={\rm div}(u_0^1\nabla u_0^2\times \nabla u_0^3)$$
yields the alternate formula
$$\int_{\Omega(0)}\det\nabla u_0(x)\,dx=\int_{\partial\Omega(0)}u_0^1\det(\nu,\nabla u_0^2,\nabla u_0^3)\,ds(x),$$
which shows that $J$ depends on $u_0^1$ only through its restriction  to the boundary. By symmetry, the same is true for every component $u_0^j$.

\bigskip

To prove Theorem \ref{th:unilat}, we begin with a property which must be rather classical. A simpler version, can be found in Milnor's paper \cite{Mil}. Notice that this lemma is valid for every parameter $\gamma>1$.
\begin{lemma}
\label{l:poly}
As a function of time, the volume $\Omega(t)$ is a polynomial of degree $\le d$.
\end{lemma}

\bigskip

\bepr

A function $G(t)$ is polynomial if and only if it is locally a polynomial. It is therefore enough to prove that $G$ is a polynomial of degree $\le d$ over some non-trivial interval $(-\epsilon,\epsilon)$.

Let us consider the map $\psi_t(x)=x+tu_0(x)$. Whenever $|t|\|u_0\|_{1,\infty}$ is less than $1$, $\psi_t$ is a diffeomorphism. Therefore $\psi_t(\Omega(0))$ is the bounded open domain whose boundary is $\psi_t(\partial\Omega(0))$. Because $\psi_t$ coincides with the flow of the fluid over $\partial\Omega(0)$, we find that $\psi_t(\Omega(0))$ is the bounded open domain whose boundary is $\partial\Omega(t)$, and we conclude
$$\Omega(t)=\psi_t(\Omega(0)).$$
We infer the formula
$$|\Omega(t)|=\int_{\Omega(0)}\det(\nabla\psi_t(x))\,dx=\int_{\Omega(0)}\det(I_d+t\nabla u_0)\,dx,$$
where the integrand is a polynomial in $t$ of degree at most $d$.

\enpr

\bigskip

Let us point out that the integral $J$ is precisely the coefficient of the monomial $t^d$ in this polynomial.
We now prove our theorem.

\bigskip

\bepr

Let $Q(t)$ be the polynomial $t\mapsto|\Omega(t)|$, $Q(t)=Jt^d+{\rm l.o.t.}$ Our assumption amounts to 
\begin{equation}\label{eq:small}
|\Omega(t)|\le C\left((1+|t|)^{d-1}\right),\qquad\forall t>0.
\end{equation}
Notice that if $J<0$, we actually have $Q(t)<0$ for large enough time, whence an immediate contradiction. However, it is important for the sequel to treat also the case $J=0$.

Let us recall the following inequality, which can be found in \cite{Chem,SerAIF}. It is a direct consequence of the Euler system (\ref{eq:cmass},\ref{eq:cmom}) and of the energy inequality (\ref{eq:ener})~:
\begin{equation}\label{eq:Chem}
\frac d{dt}\int_{\R^d}\left(\frac12\rho|tu-x|^2+t^2\rho e\right)\,dx\le2(1-d\kappa)t\int_{\R^d}\rho e\,dx.
\end{equation}
Remark that the quantity in the left integral is a linear combination of the density, the momentum and the mechanical energy with polynomial coefficients~:
$$\frac12\rho|tu-x|^2+t^2\rho e=\frac{|x|^2}2\rho+x\cdot(\rho u)+t^2\left(\frac12\rho|u|^2+\rho e\right)\,.$$

When $\gamma\in(1,\gamma_d]$, the factor $1-d\kappa$ is non-negative. With the Gronwall inequality, we deduce the estimate
\begin{equation}\label{eq:rhoga}
\int_{\R^d}\rho^\gamma dx={\rm cst}\cdot\int_{\R^d}\rho e\,dx=O\left((1+t)^{-2d\kappa}\right).
\end{equation}
Now, we can estimate the total mass with the help of H\"older Inequality:
\begin{eqnarray*}
M=\int_{\Omega(t)}\rho\,dx & \le & \left(\int_{\R^d}\rho^\gamma dx\right)^{1/\gamma}|\Omega(t)|^{1-1/\gamma} \\
& & = O\left((1+t)^{-2\frac{d\kappa}\gamma+(d-1)(1-\frac1\gamma)}\right)=O\left((1+t)^{\frac1\gamma-1}\right).
\end{eqnarray*}
Letting $t\rightarrow+\infty$, we deduce that $M=0$. In other words, there is no fluid at all.

\bigskip

If instead $\gamma>\gamma_d$, then $1-d\kappa<0$ and we  have
$$4\int_0^{+\infty}t\,dt\int_{\R^d}\rho e\,dx\le\frac1{d\kappa-1}\int_{\R^d}|x|^2\rho_0(x)\,dx,$$
wich gives 
$$\int_0^{+\infty}t\,dt\int_{\R^d}\rho^\gamma dx<\infty.$$
Using the same H\"older inequality as above, we infer
$$M^\gamma\int_0^{+\infty} t(1+t)^{(d-1)(1-\gamma)}dt<\infty,$$
which implies either $M=0$ (no fluid at all) or 
$$1+(d-1)(1-\gamma)<-1,$$
that is $\gamma>1+\frac2{d-1}$\,.

\bigskip

\hfill Q.E.D.

\bigskip

Remark that Theorem \ref{th:unilat} doesn't fully use that $(c,u)$ is a classical solution. We only need on the one hand that $(\rho,\rho u)$ be a weak solution of the Euler system, satisfying the ``entropy'' inequality (\ref{eq:ener}) (this is admissibility in the sense of Lax), and on the other hand that $(c,u)$ is of class ${\cal C}^1$ in a neighborhood of the front $\Gamma$. There is a statement analogous to Theorem \ref{th:unilat} in this context~:
\paragraph{Theorem \protect\ref{th:unilat}'}
{\em We assume $\gamma\le1+\frac1{d-1}\,$.

Let $(\rho,\rho u)$ be a  locally bounded admissible (in the sense that (\ref{eq:ener})) solution of the Euler system (\ref{eq:cmass},\ref{eq:cmom}) over $(0,+\infty)\times\R^d$. Assume that $\rho$ is supported in $\Omega$, a domain bounded in space at each time, whose boundary is smooth. Assume finally that $(c,u)$ is of class ${\cal C}^1$ in a neighborhood of $\partial\Omega$. Then necessarily}
\begin{equation}\label{eq:jacnec}
\int_{\Omega(0)}\det\nabla u_0(x)\,dx>0.
\end{equation}

 This shows that the development of shock waves is not sufficient to resolve the lack of smooth solutions. Many flows actually exhibit also some boundary singularity after some time. This is where T.-P. Liu's notion of {\em physical singularity} comes into play.

\subsection*{The mono-atomic case}

If $\gamma=\gamma_d$, Theorem \ref{th:unilat} can be improved in quantitative way. Suppose that $J\ge0$.  From
$$\frac d{dt}\int_{\R^d}\left(\frac12\rho|tu-x|^2+t^2\rho e\right)\,dx\le0,$$
we infer
$$A\int_{\Omega(t)}\rho^\gamma dx=\int_{\Omega(t)}(\gamma-1)\rho e\,dx\le\frac{\gamma-1}{t^2}\,I,\qquad I:=\int_{\Omega(0)}\frac12\rho_0|x|^2dx.$$
Applying again the H\"older Inequality, this yields
$$\frac{AM^\gamma}{\gamma-1}\,t^2\le|\Omega(t)|^{2/d}I.$$
With $|\Omega(t)|=Jt^d+{\rm l.o.t.}$, we deduce that for a solution to be smooth at the boundary for all $t>0$, we must have
$$\frac{AM^\gamma}{\gamma-1}\le J^{2/d}I.$$
We therefore have
\begin{thm}\label{th:mono}
We assume $\gamma=1+\frac1{d-1}\,$.

Let $(c_0,u_0)\in H^s(\R^d)$ be an initial data, with  $s>1+\frac d2\,$. Assume that the domain 
$$\Omega(0)=\{c_0>0\}$$ 
is bounded, with smooth boundary. Assume that 
\begin{equation}\label{eq:jacmono}
J:=\int_{\Omega(0)}\det\nabla u_0(x)\,dx<\left(\frac{2AM^\gamma}{(\gamma-1)I}\right)^{d/2}.
\end{equation}
Then the maximal time $T^*>0$ of existence provided by Theorem \ref{th:smooth} is finite:
$$T^*<+\infty.$$
\end{thm}

\subsection*{The odd-dimensional case}

As anounced above, we prove that an obstruction occurs in the much more general context of an odd space dimension, a case which includes the realistic $3$-dimensional space. 
The following result answers a question raised in \cite{Ser_five} about eternal solutions.
\begin{thm}[$d$ odd.]\label{th:nonex}
We assume that the space dimension $d$ is odd, and that $\gamma\le1+\frac1{d-1}\,$.

Then there does not exist  a non-trivial eternal solution $(c,u)\in{\cal C}^1(\R^{1+d})$ of (\ref{eq:qlmass},\ref{eq:qlmom}) (hence smooth up to vacuum), with $c_0$ supported in a  compact domain with smooth boundary.
\end{thm}

\bigskip

\bepr

Suppose $(c,u)$ is such a solution, and denote $Q(t)=|\Omega(t)|$. From Lemma \ref{l:poly}, $Q$ is a polynomial of degree $\le d$.
Because $Q$ takes positive values for every $t\in\R$, we find that  $\deg Q$ must be an even number. Because $d$ is odd, we deduce that actually $\deg Q\le d-1$, whence
\begin{equation}\label{eq:deg}
|\Omega(t)|=O\left((1+|t|)^{d-1}\right).
\end{equation}

The rest of the proof is exactly the same as in the proof of Theorem \ref{th:unilat}.

\bigskip

\hfill Q.E.D.

\bigskip

Again, Theorem \ref{th:nonex} doesn't fully use that $(c,u)$ is a classical solution, but only that it is smooth up to the vacuum. It tells us more about the onset of physical singularity at the boundary, which we consider in the next Section.

\section{Physical singularity}\label{s:phy}

It was first recognized by T.-P. Liu \cite{Liu_JJ} that many flows that are classical solutions in the interior of $\Omega(t)$ must experience a singularity at the boundary. His first motivation was the fact that for the Euler system with damping, the front travels so slowly that $\Omega(t)$ remains uniformly bounded. For gaseous stars, the gravity has a similar effect \cite{Mak}. This lack of dispersion comes in conflict with the decay of the integral of $\rho^\gamma$ and the conservation of the total mass. We have shown in  Theorem \ref{th:unilat}' that such an obstruction is present even without damping. Since the gas must flow anyway, the way to resolve the obstruction is too admit that $c$ is only H\"older-regular, of exponent $\frac12\,$, at the boundary. In other words, $c^2$ is Lipschitz and the normal derivative does not necessarily vanish
\begin{equation}\label{eq:normal}
g:=-\,\frac1{\gamma-1}\,\frac{\partial c^2}{\partial\nu}\ge0.
\end{equation}
That $g$ is non-negative follows from the fact that $c^2$ is positive in $\Omega(t)$ and vanishes on the boundary.
Assuming that $u$ is smooth accross the boundary, we obtain the identity
$$\frac{d}{dt}u(X(t),t)=g\nu$$
along a boundary path. Therefore $g$ can be viewed as an acceleration of the front expanding in vaccum. The possibility that $g$ be non-zero resolves of course the obstruction raised by Theorem \ref{th:unilat}'. It allows the volume of the domain to grow fast enough, at least as $t^d$ whatever the initial velocity. Notice that we don't need that $g$ be positive for every time. For instance, an appropriate initial data in even dimension yields a global smooth solution (Theorem \ref{th:smooth}), for which we just have $g\equiv0$.

In one space-dimension, the local existence of a flow with physical singularity at vacuum was proved in \cite{LY2} and \cite{CS1} by two different methods. These works deal with the regime where $g$ is strictly positive~; it seems to be an open problem to have an existence result covering both regimes $g=0$ and $g>0$, and in particular the transition from one regime to the other. The local existence in several space dimensions is proved in \cite{JM,CS2}.

Global existence of a weak entropy solution in one space dimension, in presence of vacuum, has been proved by means of Compensated Compactness, see \cite{DiP,DCL}. However, the method of proof says nothing about the behaviour of the solution at the front with vacuum. At least, this approach has the merit to deal with both regimes, the accelarated and the non-accelerated ones. One drawback is that the flow is obtained as the limit of approximated solutions, whose support grows faster than that expected for the genuine flow.

\subsection{One-dimensional physical singularity}
We wish to give an accurate definition of what is an admissible flow, from the point of view of physical singularity, at least in the one-dimensional case. Let $\Omega(t)=(a(t),b(t))$ be the domain occupied by the fluid (where $\rho(\cdot,t)$ is positive)\footnote{We do not consider the case where the number of connected components varies with time.}. We have $b'(t)=u(b(t),t)$ and therefore $b''=g\ge0$~; hence $b$ is convex. Likewise $a$ is a concave function. Notice that the discontinuity points of $a'$ or $b'$ form sets that are either finite or denumerable.

Let us begin with the mono-atomic case, where $\gamma=3$. We have already seen that the corresponding Euler system splits into two coupled Burgers equation. Therefore, as long as the flow is smooth, the characteristics are straight lines. Even more, the tangent lines to the boundary are characteristics. The tangent at a point $P$ actually splits into two halves, one being a $1$-characteristic and the other a $2$-characteristic. This is clear when considering a point $Q$ close to the boundary, say to the right component $\Gamma_r$~: because $b$ is convex, there are exactly two tangents to $\Gamma_r$ passing through $Q$, which are the $1$- and the $2$-characteristics respectively. We observe that a $2$-characteristic cannot emanate from $\Gamma_r$, and a $1$-characteristic cannot terminate in $\Gamma_r$.

When $\gamma>1$ is a general parameter, our definition of an admissible flow is that the same pattern be true.
\begin{defin}[$d=1$.]\label{def:adm}
Let $\gamma$ be $>1$. An admissible flow surrounded by vacuum is a measurable bounded field $(\rho\ge0,u)$, which satisfies  the following requirements:
\begin{enumerate}
\item The domain $\Omega$ is bounded at left and right by two curves $t\mapsto a(t),b(t)$, with $-a$ and $b$ convex,
\item The flow is a distributional solution of (\ref{eq:cmass},\ref{eq:cmom}) in $\Omega$,
\item It satisfies the energy inequality (\ref{eq:ener}),
\item For almost every boundary point $P\in\partial\Omega$, 
$$\lim_{(x,t)\rightarrow P}\rho(x,t)=0,$$
\item If a $2$-characteristic $\beta_2$ reaches $\Gamma_r$ at some point $P$, then $P$ is the terminal point of $\beta_2$~; if $\beta_2$ reaches $\Gamma_\ell$ at some point $Q$, then $Q$ is its initial point. Symmetrically, if a $1$-characteristic $\beta_1$ reaches $\Gamma_r$ at some point $P$, then $P$ is the initial point of $\beta_1$~; if $\beta_1$ reaches $\Gamma_\ell$ at some point $Q$, then $Q$ is its terminal point. 
\end{enumerate}
\end{defin}
The reader will have noticed that the definition above is somewhat sloppy~: without some regularity,  characteristic curves may not be well-defined.
We shall therefore use it only in situations where $(\rho,u)$ are smooth enough. This will be the case  in Theorem \ref{th:size}. Notice that our definition is consistent with the rarefaction waves displayed in the ``sub/supersonic'' example below~; this is a case where infinitely many characteristics emanate from the same boundary point.

\subsection*{Examples}

We give below two explicit examples of eternal flows with physical singularity at vacuum. They are built in one-space dimension for a mono-atomic gas ($\gamma=\gamma_1=3$), for which the Euler system decouples as a pair of Burgers equations away from shock waves. Without loss of generality, we set $A=\frac13$, so that $c=\rho$. The wave velocities $\lambda_\pm=u\pm \rho$ satisfy 
\begin{equation}\label{eq:lamBu}
\partial_t\lambda+\lambda\partial_x\lambda=0.
\end{equation}
In a domain where $\lambda$ is Lipschitz, the characteristic curves, on which $\lambda$ is constant, are lines of slope $\lambda$.

We point out that in both examples, $u_0$ vanishes identically. Because our solutions are smooth in $\Omega$, this implies the reversibility~:
$$\rho(x,-t)=\rho(x,t),\qquad u(x,-t)=-u(x,t).$$

\begin{description}
\item[An accelerated case.]
In the following example, the front between the gas and vacuum is the hyperbola defined by the equation
$$x^2=1+t^2.$$
In the gas, the density and velocity of the flow are given by
$$\rho(x,t)=\frac{\sqrt{1+t^2-x^2}}{1+t^2}\,,\qquad u(x,t)=\frac{tx}{1+t^2}\,.$$
We leave the reader verifying that each of the functions
$$\lambda_\pm(x,t)=\frac{tx\pm\sqrt{1+t^2-x^2}}{1+t^2}\,$$
satisfies the Burgers equation. Obviously, 
$$c^2=\frac{1+t^2-x^2}{(1+t^2)^2}$$
vanishes at the boundary, where it is Lipschitz, and the acceleration
$$g(\sqrt{1+t^2},t)=\frac1{(1+t^2)^{3/2}}$$
is positive for all time. Mind however that this expression is integrable in time. The domain $\Omega(t)$ behaves asymptotically as the interval $(-|t|,|t|)$, and the flow is asymptotic to 
$$\rho_R(x,t)=\frac1{|t|}\sqrt{1-\frac{x^2}{t^2}}\, ,\qquad u_R(x,t)=\frac xt\,.$$
The physical singularity is damped as $|t|\rightarrow\infty$~; the acceleration $g$ decays like $|t|^{-3}$.

Let us point out that this example does not contradict Proposition \ref{p:und}, because even after extending $c$ by zero outside of $\Omega$, the regularity of the solution $(c,u)$ is too low to implement Lax' calculation. As a matter of fact, the Cauchy--Lipchitz theorem does not apply to the characteristic flow because the wave velocities are not Lipschitz~: $\Gamma$ is the envelop of the characteristic lines.
\item[A discontinuous example.] 
This one is mimics a situation described in \cite{GB}. Let us choose instead the initial data 
$$\rho_0\equiv1,\quad u_0\equiv0\quad\hbox{in }\Omega(0)=(-1,1).$$
The solution $(c,u)$ of (\ref{eq:qlmass},\ref{eq:qlmom}) is Lipschitz, except at the points $(x,t)=(\pm1,0)$ where a discontinuity happens (obvious in the data above).  At time $t>0$, the gas occupies the domain $\Omega(t)=(-t-1,t+1)$ and the solution is given by
$$u+\rho=\left\{\begin{array}{lcr} \frac{x+1}t & \hbox{if} & -t-1<x<t-1, \\ 1 & \hbox{if} & t-1\le x<t+1. \end{array}\right.
\qquad
u-\rho=\left\{\begin{array}{lcr} -1 & \hbox{if} & -t-1<x\le-t+1, \\ \frac{x-1}t & \hbox{if} & -t+1<x<t+1. \end{array}\right.$$
The flow is a rarefaction wave in the domain defined by $|t-1|<x<t+1$, as well as in that defined by $-t-1<x<-|t-1|$. The front is accelerated only at time $t=0$, but this acceleration is an impulse~: the velocity of the front flips from $-1$ to $+1$.
\end{description}

\subsection{The domain in one space dimension}

Let us consider an eternal  flow surrounded by vacuum in the sense of definition \ref{def:adm}. We suppose in addition that $\rho ,u$ are smooth in $\Omega$. We are interested in the growth of $|\Omega(t)|$ as $t\rightarrow+\infty$. Because of the decay of the integral of $\rho^\gamma$, we already know that $|\Omega(t)|$ is bounded below by $C_1t^d$ for some constant $C_1>0$. 

\bigskip

For a one-dimensional gas, the domain $\Omega(t)=(a(t),b(t))$ is an interval where $t\mapsto -a,b$ are convex functions. The derivative $b'$ has limits $q_-\le q_+$ as $t\rightarrow\pm\infty$. Likewise $a'$ has limits $p_+\le p_-$. Because of the lower bound of the volume, we know that
$$p_+<q_+,\qquad q_-<p_-.$$
We actually have
\begin{thm}[$d=1$.]\label{th:size} Let a flow be smooth (shock-free) in its domain $\Omega(t)$, where $\Omega(0)=(a_0,b_0)$ is a bounded interval. We assume that the flow is admissible in the sense of Definition \ref{def:adm}.

Then $p_+=q_-$ and $p_-=q_+$. 
\end{thm}

\bepr

We consider $2$-characteristics, which are integral curves of the ODE
$$\frac{dX}{dt}=\lambda(X,t),\qquad\lambda:=u+c.$$
Because the flow is smooth in $\Omega$, we know that $r:=u+\frac c\kappa$ is constant along each $2$-characteristic.

We choose the origin of the time arrow in such a way that $a'$ and $b'$ be continuous at $t=0$.
Denote $\Omega^+=\Omega\cap\{t>0\}$. By Cauchy--Lipschitz, $\Omega^+$ is foliated by the $2$-characteristics. On such a curve $\beta$, the time goes from $T_{in}\ge0$ to $T^{fin}$. At $t=T_{in}$, $\beta$ reaches the boundary $\partial\Omega^+$ at some point $m$. By admissibility, we actually have
$$m\in\Gamma^+:=\Gamma_\ell^+\cup(a_0,b_0),\qquad\Gamma_\ell^+:=\Gamma_\ell\cap\{t\ge0\}.$$
We order the set of $2$-characteristics in $\Omega^+$ from left to right, and denote it $(-\infty,\beta_0)$, where $\beta_0=\{(b_0,0)\}$ and $-\infty$ is the ``point at infinity'' along $\Gamma_\ell^+$.

Consider two such characteristics $\beta\prec\beta'$. If $T^{fin}(\beta)$ is finite, the terminal point $\beta(T^{fin})$ is on $\Gamma_r$. Then $\beta'$ is contained in the bounded set delimited by $\beta$ at left and $\Gamma_r$ at right, and we infer $T^{fin}(\beta')\le T^{fin}(\beta)<\infty$. We deduce that $\beta\mapsto T^{fin}$ is non-increasing, and there exists  $\beta^*\in[-\infty,\beta_0]$ such that 
\begin{eqnarray*}
(\beta\prec\beta^*) & \Longrightarrow & (T^{fin}(\beta)=+\infty), \\
(\beta^*\prec\beta) & \Longrightarrow & (T^{fin}(\beta)<+\infty).
\end{eqnarray*}
In the limit case where $T^{fin}$ is finite for every $\beta$ (unlikely), then $\beta^*=-\infty$. If instead  $T^{fin}=+\infty$ for every $\beta$ (unlikely), then $\beta^*=b_0$.

When $\beta\prec\beta^*$, the Ricatti equation along $\beta$, plus the fact that
$$\frac{\partial\lambda}{\partial r}=\frac{\gamma+1}4>0,$$
imply that $\partial_xr\ge0$. This shows that $\beta\mapsto r|_\beta$ is non-decreasing up to $\beta^*$. On the other hand, if $\beta^*\prec\beta$, we have seen that $\beta\mapsto T^{fin}(\beta)$ is non-increasing. Since $r$ equals $u=b'$ over the boundary, and $b'$ is non-decreasing, we deduce that $\beta\mapsto r|_\beta$ is non-increasing beyond $\beta^*$. If $\beta^*\prec b_0$, we have
$$q_+=\sup_{\Gamma^+_r}u\le\sup_{\Gamma^+_r}r\le\sup_{\Omega^+}r=r|_{\beta^*}
=\sup_{\beta^*\prec\beta}r|_\beta\le\sup_{\Gamma_r^+}u= q_+.$$
If instead $\beta^*=b_0$, then 
$$q_+=\sup_{\Gamma^+_r}u\le\sup_{\Gamma^+_r}r\le\sup_{\Omega^+}r=\lim_{\beta\rightarrow \beta_0}r|_\beta=\lim_{x\rightarrow b_0}r(x,0)=r(b_0,0)=b'(0)\le q_+,$$
where the last equality stands because $b'$ is continuous at $t=0$.
In all cases, we infer
$$\sup_{\Omega^+}r=r|_{\beta^*}=q_+.$$

Considering symmetrically the foliation of $\Omega_-=\Omega\cap\{t<0\}$ by the $2$-characteristics, we find a $\beta_*\in[\alpha_0,+\infty]$ such that
\begin{eqnarray*}
(\beta\prec\beta_*) & \Longrightarrow & (T^{in}(\beta)>-\infty), \\
(\beta_*\prec\beta) & \Longrightarrow & (T^{in}(\beta)=-\infty).
\end{eqnarray*}
We have that $\beta\mapsto r|_\beta$ is non-decreasing over $(a_0,\beta_*)$, non-increasing over $(\beta_*,+\infty)$, and satisfies
$$\sup_{\Omega^-}r=r|_{\beta_*}=p_-.$$

We point out that those $2$-characteristics passing through a point $m\in(\alpha_0,b_0)$ belong to both sets, associated with either $\Omega^+$ or $\Omega_-$. 

Suppose now that $p_-\ne q_+$. Without loss of generality, we have $p_-<q_+$. Then $r|_{\beta^*}>r|_\beta$ for every $2$-characteristic $\beta$ in $\Omega_-$. Thus $\beta^*$ does not emanate from a point $m\in(a_0,b_0)$~; we must have $\beta^*\prec \alpha_0$. Thus $\beta^*$ emanates from some point $(a(t),t)$ of $\Gamma_\ell^+$ and we have $r|_{\beta^*}=r(a(t),t)$. This gives $q_+=a'(t)\le p_-$, a contradiction. 

This ends the proof of the Theorem.

\enpr

\subsection*{Open problems}

One important issue is of course to give a multi-dimensional version of Definition \ref{def:adm}. The difficulty here is that we cannot distinguish between two characteristic families~: the characteristic curves passing through a point $P$ are tangent to a cone, which is a connected set. This cone shrinks to a line at boundary points.

\bigskip

Because of the dispersion estimate and the conservation of mass, we know that $|t|^d=O(|\Omega(t)|)$ as $|t|\rightarrow+\infty$. We suspect that the converse inequality $|\Omega(t)|=O(t^d)$ is also true. Theorem \ref{th:size} provides a positive answer in the special case where $d=1$ and the flow is shock-free. It is also true for flows that are smooth up to the boundary, see Lemma \ref{l:poly}. 

Actually, in both cases, the domain $\Omega(t)$ behaves asymptotically as $t\cal O$ for the {\em same} $\cal O$ as $t\rightarrow\pm\infty$. This latter property is doubtful for flows that display shocks, but is it true for shock-free flows with physical singularity when $d\ge2$~?

\paragraph{Acknowledgments.} I thank warmly Alexis Vasseur who welcame me at the Department of Mathematics of the University of Texas at Austin~; I came back to the topic of  eternal solutions when I learnt from him about Levermore's construction of global Maxwellians. I am grateful to \'Etienne Ghys and Bruno S\'evennec, who told me about the similarity of one of my arguments above with the elegant proof by J. Milnor of the hairy ball theorem. I am indebted to Constantin Dafermos, who refreshed my memory about Liu's work on the physical singularity.

\end{document}